\documentclass[11pt]{article}

\usepackage{amsmath}
\usepackage{amssymb}
\usepackage{eucal}

\begin{document}

\baselineskip 16pt

\title{\bf  One application of the 
$\sigma$-local formations  of finite groups }

\author{Zhang Chi \thanks{Research of the first author is supported by 
 China Scholarship Council and NNSF of 
China(11771409)}\\
{\small Department of Mathematics, University of Science and
Technology of China,}\\ {\small Hefei 230026, P. R. China}\\
{\small E-mail:
zcqxj32@mail.ustc.edu.cn}\\ \\
{ Alexander  N. Skiba}\\
{\small Department of Mathematics and Technologies of Programming,
  Francisk Skorina Gomel State University,}\\
{\small Gomel 246019, Belarus}\\
{\small E-mail: alexander.skiba49@gmail.com}}

 \date{}
\maketitle

\begin{abstract}  Throughout this paper, all groups are finite.
 Let $\sigma =\{\sigma_{i} |
 i\in I \}$  be some partition of the set of all primes 
$\Bbb{P}$.   If
 $n$ is an integer, the symbol $\sigma (n)$ denotes
 the set  $\{\sigma_{i} |\sigma_{i}\cap \pi (n)\ne 
 \emptyset  \}$. The   integers $n$ 
and $m$ are called $\sigma$-coprime if $\sigma (n)\cap \sigma (m)=\emptyset$.

   Let $t > 1$ be a natural number and  let  
$\mathfrak{F}$ be a class  of groups.  Then we say that $\mathfrak{F}$ is  
$\Sigma_{t}^{\sigma}$-closed provided  $\mathfrak{F}$ contains each group
 $G$  with  subgroups $A_{1},  \ldots , A_{t}\in
 \mathfrak{F}$ whose  indices
$|G:A_{1}|$,   $\ldots$,  $|G:A_{t}|$ are pairwise  $\sigma$-coprime.

In this paper,    we study $\Sigma_{t}^{\sigma}$-closed classes of finite  groups.

\end{abstract}


\let\thefootnote\thefootnoteorig

\section{Introduction}

Throughout this paper, all groups are finite and $G$ always denotes
a finite group. Moreover,  $\mathbb{P}$ is the set of all  primes,
  $\pi= \{p_{1}, \ldots , p_{n}\} \subseteq  \Bbb{P}$ and  $\pi' =  \Bbb{P} \setminus \pi$. If
 $n$ is an integer, the symbol $\pi (n)$ denotes
 the set of all primes dividing $n$; as usual,  $\pi (G)=\pi (|G|)$, the set of all
  primes dividing the order of $G$.

Following Shemetkov \cite{26}, $\sigma$  is some partition of  
$\Bbb{P}$, that is,  $\sigma =\{\sigma_{i} |
 i\in I \}$, where   $\Bbb{P}=\bigcup_{i\in I} \sigma_{i}$
 and $\sigma_{i}\cap \sigma_{j}= \emptyset  $ for
 all $i\ne j$; $\Pi \subseteq \sigma$ and $\Pi' = \sigma \setminus \Pi
$.   The group  $G$ is said to be  \cite{1}:   
\emph{$\sigma$-primary}  if  $G$ is a $\sigma_{i}$-group for some 
$i$; \emph{$\sigma$-soluble} 
 if every chief factor of $G$ is $\sigma$-primary.

In what follows,     $\sigma (n)=\{\sigma_{i} |\sigma_{i}\cap \pi (n)\ne 
 \emptyset  \}$  \cite{1, alg12},    $\sigma (G)=\sigma (|G|)$    and $\sigma
 (\mathfrak{F})=\bigcup_{G\in \mathfrak{F}}\sigma (G)$. The  the integers $n$ 
and $m$ are called \emph{$\sigma$-coprime} if $\sigma (n)\cap \sigma (m)=\emptyset$.

Recall also that $G$ is called \emph{$\sigma$-decomposable}  (Shemetkov \cite{26})  or 
\emph{$\sigma$-nilpotent}   (Guo and Skiba \cite{33}) if $G=G_{1}\times \cdots \times G_{n}$ 
for some $\sigma$-primary groups $G_{1}, \ldots, G_{n}$;   
\emph{meta-$\sigma$-nilpotent}   \cite{comm} if $G$ is an extension of 
some $\sigma$-nilpotent group be the $\sigma$-nilpotent group.

The $\sigma$-nilpotent groups  have   proved to be very useful in
 the formation theory (see, in particular, the papers    \cite{19, 
20} and  the books \cite[Ch. IV]{26}, \cite[Ch. 6]{BalE}). 
  In the recent years, the $\sigma$-nilpotent groups and 
 various classes of meta-$\sigma$-nilpotent groups  have found new and to
 some extent unexpected  
applications in the theories of permutable and generalized subnormal 
subgroups  (see, in particular,  \cite{1, alg12},    \cite{3}--\cite{6} and the survey \cite{comm}).
This circumstance make the task of further studying   of  $\sigma$-nilpotent and 
  meta-$\sigma$-nilpotent groups   quite actual and interesting.

In this paper, we study $\Sigma_{t}^{\sigma}$-closed classes of meta-$\sigma$-nilpotent groups
 in the sense of the following

 { \bf Definition 1.1.}  Let $t > 1$ be a natural number and  let  
$\mathfrak{F}$ be a class  of groups.  Then we say that $\mathfrak{F}$ is  
\emph{$\Sigma_{t}^{\sigma}$-closed} provided  $\mathfrak{F}$ contains each group
 $G$  with  subgroups $A_{1},  \ldots , A_{t}\in
 \mathfrak{F}$ whose  indices
$|G:A_{1}|$,   $\ldots$,  $|G:A_{t}|$ are pairwise  $\sigma$-coprime.

We will omit the symbol $\sigma$ 
in this definition in the classical case, 
 when $\sigma =\sigma^{1} = \{\{2\}, \{3\}, \ldots \}$  (we use here the notation in 
\cite{alg12}).  Thus in this case we deal with  
\emph{$\Sigma_{t}$-closed} classes of groups, in the usual sense  
(see L.A. Shemetkov \cite[p. 44]{26}).

Recall that a class of groups  $\mathfrak{F}$ is called a \emph{formation} if:
 (i) $G/N\in \mathfrak{F}$ whenever  $G\in \mathfrak{F}$, and (ii)
 $G/(N\cap R)\in \mathfrak{F}$ whenever
  $G/N\in \mathfrak{F}$ and  $G/R\in \mathfrak{F}$.   
The formation  $\mathfrak{F}$
 is called   \emph{saturated} or \emph{local} if 
 $G\in \mathfrak{F}$ whenever $G/\Phi (G)\in \mathfrak{F}$.

We call  any function $f$ of the form
$$f:\sigma  \to\{\text{formations of groups}\}$$
\emph{a formation $\sigma$-function} \cite{problem4}, and we put $$LF_{\sigma}(f)=(G
 \mid  G=1  \  \text{or   }\ G\ne 1\ \text{and   }\  G/O_{\sigma _{i}',   \sigma _{i}}(G) \in
 f(\sigma 
_{i})  \  \text{for all  } \sigma _{i} \in  \sigma(G)).$$

{\bf Definition 1.2.}   If for some
 formation $\sigma$-function $f$ 
 we have $\mathfrak{F}=LF_{\sigma}(f)$, then we say, following \cite{problem4}, that the class  
$\mathfrak{F}$ is \emph{ $\sigma $-local } and $f$ is a
 \emph{ $\sigma$-local definition of
$ \mathfrak{F}$.}

Before continuing, consider some examples.

{\bf Example 1.3.}  (i) 
In view of \cite[IV, 3.2]{DH},
 in the  case 
 when $\sigma =\sigma^{1}$, a  formation $\sigma$-function and a  $\sigma $-local formation are, respectively, 
 a formation function and a local formation  in the usual sense 
\cite[IV, Definition 3.1]{DH} (see also \cite[Chapter 2]{BalE}).
 We use in this case  instead of $LF_{\sigma}(f)$  the symbol  $LF(f)$, as 
usual \cite[IV, Definition 3.1]{DH}.

 (ii) For  the formation  of all identity groups ${\mathfrak{I}}$ we have
 ${\mathfrak{I}}=LF_{\sigma}(f)$, where $f(\sigma _{i})=\emptyset $ for all $i$. 

(iii) Let ${\mathfrak{N}}_{\sigma}$ be the class of all $\sigma$-nilpotent 
groups. Then ${\mathfrak{N}}_{\sigma}$  is a formation \cite{1} and, clearly,  
${\mathfrak{N}}_{\sigma}=LF_{\sigma}(f)$, where $f(\sigma _{i})={\mathfrak{I}}$ 
 for all $i$. 

(iv)  Now let ${\mathfrak{N}}_{\sigma}^{2}$ be the class of all 
meta-$\sigma$-nilpotent groups. Then  
${\mathfrak{N}}_{\sigma}^{2}=LF_{\sigma}(f)$, where $f(\sigma _{i})={\mathfrak{N}}_{\sigma}$ 
 for all $i$. 

(v) The formation of all supersoluble groups ${\mathfrak{U}}$ is not 
$\sigma$-local for every $\sigma$  with $\sigma \ne \sigma ^{1}$. Indeed, 
suppose that ${\mathfrak{U}}=LF_{\sigma}(f)$ is  
$\sigma$-local and for some $i$ we have $|\sigma _{i}| > 1$. Let $p, q\in 
\sigma _{i}$, where $ p > q$. Finally,  let $G=C_{q}\wr C_{p}=K\rtimes C_{p}$ be the 
regular wreath product of groups $ C_{q}$  and $ C_{p}$ with $ |C_{q}|=q$  and 
$ |C_{p}|=p$, where $K$ is the base group of $G$. Then $C_{G}(K)=K$ and, also,  
$O_{\sigma _{i}', \sigma _{i}}(G)=G$ and $\sigma 
(G)=\{\sigma _{i}\}$. Since $C_{p}\in {\mathfrak{U}}$,  
$f(\sigma _{i})\ne \emptyset$.  
Hence  $G\in LF_{\sigma}(f)={\mathfrak{U}}$, so $G=C_{q}\times C_{p}$ since $ p > q$, a contradiction.
Hence we have (iv).

The theory of 
$\Sigma_{t}$-closed classes of soluble groups   and various its applications
 were considered by 
Otto-Uwe Kramer in \cite{kram} (see also \cite[Chapter 1]{26} or \cite[Chapter 2]{GuoI}).

Our  main goal here is to prove the following result.

{\bf Theorem 1.4.}   {\sl Every $\sigma$-local formation of meta-$\sigma$-nilpotent groups  
is   $\Sigma_{4}^{\sigma}$-closed.}

In the case when $\sigma =\sigma^{1}$, we get from Theorem 1.4 the 
following well-known facts.

{\bf Corollary 1.5} (Doerk \cite{Dok1}).  {\sl If $G$ has four supersoluble 
 subgroups $A_{1},  A_{2}, A_{3}, A_{4}$  whose   indices 
$|G:A_{1}|,  |G:A_{2}|, |G:A_{3}|,  |G:A_{4}|$
are pairwise  coprime, then  $G$ is itself  supersoluble. }

{\bf Corollary 1.6.}  {\sl If $G$ has four meta-nilpotent 
 subgroups $A_{1},  A_{2}, A_{3}, A_{4}$  whose   indices 
$|G:A_{1}|,  |G:A_{2}|, |G:A_{3}|,  |G:A_{4}|$
are pairwise  coprime, then  $G$ is itself  meta-nilpotent. }

{\bf Corollary 1.7.} {\sl Suppose that  $G$ has four 
 subgroups $A_{1},  A_{2}, A_{3}, A_{4}$  whose   indices 
$|G:A_{1}|,  |G:A_{2}|, |G:A_{3}|,  |G:A_{4}|$
are pairwise  coprime. If  the derived subgroup $A_{i}'$ of  $A_{i}$ is
  nilpotent  for all $i=1, 2, 3, 4 $, then $G'$ is nilpotent. }

Finally, we get from Theorem 1.4 the following

 {\bf Corollary 1.8} (Otto-Uwe Kramer \cite{kram}).
 {\sl  Every local formation of
 meta-nilpotent groups is $\Sigma_{4}$-closed.}

In fact, in the theory of the $\pi$-soluble groups ($\pi =\{p_{1}, \ldots , p_{n}\}$)
  we deal with  the 
partition 
 $\sigma =\sigma ^{1\pi }=\{\{p_{1}\},  \ldots , \{p_{n}\}, \pi' \}$ of $\Bbb{P}$ \cite{alg12}.  
 Note that    $G$ is:     $\sigma ^{1\pi }$-soluble
 if and only if $G$ is    $\pi$-soluble;       $\sigma ^{1\pi }$-nilpotent
 if and only if $G$ is    \emph{$\pi$-special} \cite{Cun2},
 that is, $G=O_{p_{1}}(G)\times \cdots \times
 O_{p_{n}}(G)\times O_{\pi'}(G)$.  
Hence in this case we get from  Theorem 1.4 the following  results.

{\bf Corollary 1.9.}   {\sl  Suppose that  $G$ has four meta-$\pi$-special
 subgroups $A_{1},  A_{2}, A_{3}, A_{4}$  whose   indices 
$|G:A_{1}|,  |G:A_{2}|, |G:A_{3}|,  |G:A_{4}|$
are pairwise  coprime and   each of them is either a $\pi$-number or a $\pi'$-number.
 Suppose also that at most  one of the numbers
 $|G:A_{1}|, |G:A_{2}|, |G:A_{3}|,    |G:A_{4}|$
is a $\pi'$-number.  Then $G$ is  meta-$\pi$-special. }

{\bf Corollary 1.10.}   {\sl Suppose that $G$ has   subgroups $A_{1},  \ldots ,
 A_{4}$   such that   the  indices 
$|G:A_{1}|,  \ldots ,  |G:A_{4}|$
are pairwise  coprime and   each of them is either a $\pi$-number or a $\pi'$-number. 
  Suppose also that at most  one of the numbers
 $|G:A_{1}|,  |G:A_{2}|, |G:A_{3}|,    |G:A_{4}|$
is a $\pi'$-number. 
  If the derived subgroup $A_{i}'$ of  $A_{i}$ is
  $\pi$-special for all $i$, then $G'$ is  $\pi$-special. }

If for a subgroup $A$ of $G$ we have $\sigma (|A|) \subseteq \Pi$
 and $\sigma (|G:A|) \subseteq \Pi'$, then $A$ is said to be a   \emph{Hall 
$\Pi$-subgroup} \cite{comm} of $G$.  
 We say  that $G$ is \emph{ $\Pi$-closed} if $G$ has a normal Hall 
$\Pi$-subgroup.

The proof of Theorem 1.4 is preceded by a large number of auxiliary 
results. The following theorem is one of them.

{\bf Theorem  1.11.}     (i) {\sl The class of all $\sigma$-soluble
  $\Pi$-closed groups
 $\mathfrak{F}$  is $\Sigma_{3}^{\sigma}$-closed.   }

(ii) {\sl Every   formation of $\sigma$-nilpotent groups
   $\mathfrak{M}$   is $\Sigma_{3}^{\sigma}$-closed.  }

{\bf Corollary   1.12.}     (i) {\sl The classes  of all $\sigma$-soluble
   groups  and of all  $\sigma$-nilpotent groups  are   $\Sigma_{3}^{\sigma}$-closed.   }

  In the case when $\sigma =\sigma^{1}$, we get from Corollary  1.12 the 
following well-known results.

{\bf Corollary 1.13} (Wielandt  \cite[Ch. I, Theorem 3.4]{DH}). {\sl  If $G$ has three
soluble subgroups $A_{1}$, $A_{2}$ and $A_{3}$ whose indices
$|G:A_{1}|$, $|G:A_{2}|$, $|G:A_{3}|$ are pairwise  coprime, then
$G$ is itself soluble. }

{\bf Corollary 1.14} (Kegel \cite{Kegel1}). {\sl  If $G$ has three
nilpotent subgroups $A_{1}$, $A_{2}$ and $A_{3}$ whose indices
$|G:A_{1}|$, $|G:A_{2}|$, $|G:A_{3}|$ are pairwise  coprime, then
$G$ is itself nilpotent. }

 {\bf Corollary 1.15} (Doerk \cite{Dok1}). {\sl  If  $G$ has three
abelian subgroups $A_{1}$, $A_{2}$ and $A_{3}$ whose indices
$|G:A_{1}|$, $|G:A_{2}|$, $|G:A_{3}|$ are pairwise  coprime, then
$G$ is itself abelian. }

In the case when $\sigma =\sigma ^{1\pi }$, we get from  Theorem 1.11 the 
following facts.

{\bf Corollary 1.16.}   {\sl  Suppose that    $G$ has three $\pi$-soluble
 subgroups $A_{1},  A_{2}, A_{3}$  whose   indices 
$|G:A_{1}|,  |G:A_{2}|, |G:A_{3}|$
are pairwise  coprime and   each of them is either a $\pi$-number or a $\pi'$-number.
  Suppose also that at most  one of the numbers
 $|G:A_{1}|,  |G:A_{2}|,  |G:A_{3}|$
is a $\pi'$-number.    Then $G$ is  $\pi$-soluble. }

{\bf Corollary 1.17.}   {\sl  Suppose that    $G$ has three $\pi$-special
 subgroups $A_{1},  A_{2}, A_{3}$  whose   indices 
$|G:A_{1}|,  |G:A_{2}|, |G:A_{3}|$
are pairwise  coprime and   each of them is either a $\pi$-number or a $\pi'$-number.  Suppose also that at most  one of the numbers
 $|G:A_{1}|,  |G:A_{2}|,  |G:A_{3}|$
is a $\pi'$-number.    Then $G$ is  $\pi$-spesial. }

\section{General properties of $\sigma$-local formations
}

If $\mathfrak{M}$ and    $\mathfrak{H}$  are classes of groups,  then  
$\mathfrak{M}\mathfrak{H} $  is the class of groups $G$ such that for some 
normal subgroup $N$ of $G$ we have $G/N\in \mathfrak{H} $ and  $N\in 
\mathfrak{M} $. The \emph{Gasch\"{u}tz product}  $\mathfrak{M}\circ 
\mathfrak{H} $ of $\mathfrak{M}$ and    $\mathfrak{H}$ is defined as 
folows:  $G\in \mathfrak{M}\circ 
\mathfrak{H} $ if and only if $G^{\mathfrak{H}}\in \mathfrak{M}$.    The 
class 
 $\mathfrak{F}$
 is called   \emph{hereditary} in the sense of Mal'cev  \cite{Mal}  if 
 $G\in \mathfrak{F}$ whenever $G\leq A\in \mathfrak{F}$.

All statements of the following lemma  are well-known
 (see,  \cite[Chapter II]{shem-sk} or 
 \cite[Chapter IV]{DH}) and, in fact, each of them may be proved by the direct calculations.

{\bf Lemma  2.1.}   {\sl Let $\mathfrak{M}$,   $\mathfrak{H}$  and 
 $\mathfrak{F}$ be  formations.  }

(1)  {\sl  $\mathfrak{M}\circ\mathfrak{H}$
 is a formation.}

(2) {\sl If  $\mathfrak{M}$ is  hereditary, then
 $ \mathfrak{M}\mathfrak{H}=\mathfrak{M}\circ\mathfrak{H}$.}

(3)  $(\mathfrak{M}\circ\mathfrak{H}) \circ \mathfrak{F}=
 \mathfrak{M}\circ (\mathfrak{H}\circ\mathfrak{F})$.

(4) {\sl If  $\mathfrak{M}$ and $\mathfrak{H}$ are hereditary,
 then $\mathfrak{M}\mathfrak{H}$ is   hereditary.
}

(5) {\sl If  $\mathfrak{M}$ is saturated and $\pi (\mathfrak{H}) \subseteq \pi (\mathfrak{M})$,
 then $\mathfrak{M} \circ \mathfrak{H}$ is  saturated.
}

 We write ${\mathfrak{G}}_{\Pi}$ (respectively ${\mathfrak{S}}_{\Pi}$) 
 to denote  the class  of 
all $\Pi$-groups (respectively the class  of 
all $\sigma$-soluble $\Pi$-groups).
  In particular, 
 ${\mathfrak{G}}_{\sigma _{i}'}$ is the class  of 
all  $\sigma _{i}'$-groups and  $\mathfrak{G}_{\sigma _{i}}$ is the class  of 
all  $\sigma _{i}$-groups and  ${\mathfrak{S}}_{\sigma _{i}'}$ is the class  of 
all $\sigma$-soluble  $\sigma _{i}'$-groups.

We use $F_{\Pi}(G)$ to denote the product of all normal  $\Pi'$-closed 
subgroups of $G$; we write also $F_{\sigma _{i}}(G)$ instead of $F_{\{\sigma _{i}\}}(G)$.

{\bf Lemma  2.2.}    (1)   {\sl The class of all ($\sigma$-soluble)
  $\Pi$-closed groups
 $\mathfrak{F}$  is a hereditary  formation.   Moreover, } 

(2) {\sl  If 
$E$ is a normal  subgroup of $G$ and $E/E\cap \Phi (G)\in \mathfrak{F}$, then
  $E\in \mathfrak{F}$.   Hence the formation $\mathfrak{F}$  is saturated.}

(3) {\sl   If   
$A, B\in  \mathfrak{F}$ are normal  subgroups of $G$ and $G=AB$, then
  $G\in \mathfrak{F}$.  }

(4) {\sl   If 
$E$ is a subnormal  subgroup of $G$, then $F_{\Pi}(G)\cap E=F_{\Pi}(E)$. }

{\bf Proof. } (1)  It is clear that $\mathfrak{F}={\mathfrak{G}}_{\Pi}{\mathfrak{G}}_{\Pi'}$.  
Hence $\mathfrak{F}$  is a hereditary formation by Lemma  2.1(1, 2, 4).

(2) Let  $H/E\cap \Phi (G)$ be the normal
 Hall $\Pi$-subgroup of $E/E\cap \Phi (G)$.  Then $H/E\cap \Phi (G)$ is 
characteristic in  $E/E\cap \Phi (G)\trianglelefteq G/E\cap \Phi (G)$, so 
$H$ is normal in $G$.    
 Let $D=O_{\Pi'}(E\cap \Phi (G))$. Then, since $E\cap \Phi
 (G)$ is 
nilpotent,  $D$ is a Hall  $\Pi'$-subgroup of $H$.  
Hence by the Schur-Zassenhaus  
theorem, $H$ has a Hall $\Pi$-subgroup, say  $V$, and   any two Hall $\Pi$-subgroups of $H$
 are conjugated in $H$. Therefore,  
 $G=HN_{G}(V)=(VD))N_{G}(V)=N_{G}(V)$ by the Frattini argument.  Thus $V$   is 
normal in $G$. Finally, $V$ is a Hall $\Pi$-subgroup of $E$ since
 $ \sigma (|E/E\cap \Phi (G): H/E\cap \Phi (G)|)\cap \Pi=\emptyset,$  so 
  $E\in \mathfrak{F}$.  

(3) If  $V$ is a Hall  $\Pi$-subgroup of  $A$, then $V$ is characteristic 
in $A$ and so $V$ is normal in $G$. Similarly, a Hall  
$\Pi$-subgroup $W$ of  $B$ is normal in $G$.   Moreover, 
$$G/VW=AB/VW=(AVW/VW)(BVW/VW),$$ where  $$AVW/VW\simeq  A/A\cap VW=A/V(A\cap 
W)\simeq  (A/V)/(V(A\cap W)/V$$ and  $BVW/VW$ are $\Pi'$-groups.
 Hence $VW$ is  a Hall  $\Pi$-subgroup of  $G$, so $G\in \mathfrak{F}$.

(4)  Since  the group $A$ is $\Pi'$-closed if and only if 
$A\in {\mathfrak{G}}_{\Pi'}{\mathfrak{G}}_{\Pi}$, we have (4) by
 \cite[VIII, Proposition 2.4(d)]{DH}. 

 The 
lemma is proved.

If $f$ is a formation $\sigma$-function, then  the symbol $\text{Supp}(f)$  denotes the
  \emph{support }  of $f$, that is, the set of all 
$\sigma _{i}$  such that $f(\sigma _{i})\ne \emptyset$.

{\bf Lemma 2.3. }  {\sl Let $\mathfrak{F}=LF_{\sigma}(f)$  and $\Pi =\text{Supp}(f)$}.

 (1)      {\sl    $\Pi=\sigma (\mathfrak{F})$. }

 (2)  {\sl  $G\in \mathfrak{F}$ if and only if
 $G\in {\mathfrak{G}}_{\sigma_{i}'}{\mathfrak{G}}_{\sigma _{i}}f(\sigma _{i})$
      for all  $\sigma _{i}\in \sigma (G)$. }
 
(3)  {\sl
 $\mathfrak{F}= (\bigcap _{\sigma _{i}\in 
\Pi}\mathfrak{G}_{\sigma_{i}'} \mathfrak{G}_{\sigma _{i}}f(\sigma 
_{i})) \cap  \mathfrak{G}_{\Pi }.$ Hence $\mathfrak{F}$ is a saturated formation.}

(4) {\sl If every group in $\mathfrak{F}$ is $\sigma$-soluble, then $\mathfrak{F}= (\bigcap _{\sigma _{i}\in 
\Pi}\mathfrak{S}_{\sigma_{i}'} \mathfrak{G}_{\sigma _{i}}f(\sigma 
_{i})) \cap  \mathfrak{S}_{\Pi }.$}

{\bf Proof.} (1)  If    $\sigma _{i}\in \Pi$, then $1\in f(\sigma _{i})$ and for every 
$\sigma _{i}$-group $G\ne 1$ we have $\sigma (G)=\{\sigma _{i}\}$ and 
$O_{\sigma _{i}',   \sigma _{i}}(G)=G$. Hence $G\in LF_{\sigma}(f)=\mathfrak{F},$ so
  $\sigma _{i}\in \sigma (\mathfrak{F})$.   Therefore $\Pi\subseteq \sigma 
(\mathfrak{F})$.  On the other hand,  if  $\sigma _{i}\in \sigma (\mathfrak{F})$, then  
for some group $G\in \mathfrak{F}$ we have  $\sigma _{i}\in\sigma (G)$ and 
 $G/F_{\sigma _{i}}(G)\in f(\sigma _{i})$. Thus  $\sigma _{i}\in \Pi$, so 
  $\Pi=\sigma (\mathfrak{F})$.

(2)  If  $G\in \mathfrak{F}$ and  $\sigma _{i}\in \sigma (G)$, then 
$G/F_{\sigma _{i}}(G)\in f(\sigma _{i})$, where $ F_{\sigma _{i}}(G)$  
 is $\sigma _{i}'$-closed by Lemma 2.2(3). 
 Hence $G\in  \mathfrak{G}_{\sigma_{i}'} \mathfrak{G}_{\sigma _{i}}f(\sigma 
_{i})$ by Lemma 2.2(1).  
Similarly, if for any  $\sigma _{i}\in \sigma (G)$ we have  
 $G\in  \mathfrak{G}_{\sigma_{i}'} \mathfrak{G}_{\sigma _{i}}f(\sigma 
_{i})$, then $G/F_{\sigma _{i}}(G)\in f(\sigma 
_{i})$ and so $G\in \mathfrak{F}$. 

(3)    If $G\in \mathfrak{F}$, then $\sigma (G)\subseteq  \Pi$ and so  
 $G\in \mathfrak{G}_{\Pi }.$ Moreover, in this case for every  $\sigma _{i}\in \sigma (G)$ we have  $G\in 
{\mathfrak{G}}_{\sigma_{i}'}{\mathfrak{G}}_{\sigma _{i}}f(\sigma _{i})$
 by Part (2).  Finally, if 
$\sigma _{i}\in \Pi\setminus \sigma (G)$, then $G\in {\mathfrak{G}}_{\sigma_{i}'}\subseteq 
{\mathfrak{G}}_{\sigma_{i}'}{\mathfrak{G}}_{\sigma _{i}}f(\sigma _{i})$ 
since the class    ${\mathfrak{G}}_{\sigma_{i}'}$ is  hereditary. 
Therefore    $\mathfrak{F}\subseteq  (\bigcap _{\sigma _{i}\in 
\Pi}\mathfrak{G}_{\sigma_{i}'} \mathfrak{G}_{\sigma _{i}}f(\sigma 
_{i})) \cap  \mathfrak{G}_{\Pi }.$  Hence $\mathfrak{F}=  (\bigcap _{\sigma _{i}\in 
\Pi}\mathfrak{G}_{\sigma_{i}'} \mathfrak{G}_{\sigma _{i}}f(\sigma 
_{i})) \cap  \mathfrak{G}_{\Pi }$  is a saturated 
formation by Lemmas 2.1(5)  and 2.2(1, 2).      Hence we have (3).

(4) See the proof of (3).

The lemma is proved.

{\bf Lemma 2.4.}  {\sl If   $\mathfrak{F}=LF_{\sigma}(f)$, then 
 $\mathfrak{F}=LF_{\sigma}(t)$, where $t(\sigma _{i})=f(\sigma 
_{i})\cap \mathfrak{F}$ for all $\sigma _{i}\in \sigma .$   }

{\bf Proof.  }   First note that in view of Lemma 2.3(3),
 $t$ is a formation  
$\sigma$-function   and $LF_{\sigma}(t) \subseteq 
\mathfrak{F}$. On the other hand, if  $G\in \mathfrak{F}$, then
$G/F_{\sigma _{i}}(G)\in f(\sigma _{i})\cap \mathfrak{F}=t(\sigma _{i})$ for every $\sigma _{i}\in \sigma (G)$  
 and so $G\in LF_{\sigma}(t)$.
 Hence $\mathfrak{F}=LF_{\sigma}(t)$.    
The lemma is proved.

{\bf Proposition  2.5.} {\sl Let $f$ and $h$  be  formation $\sigma$-functions
 and let }  $\Pi=\text{Supp}(f)$. {\sl Suppose that
  $\mathfrak{F}=LF_{\sigma}(f)=LF_{\sigma}(h)$. }
 
(1) {\sl If $\sigma _{i}\in \Pi $, then
  $\mathfrak{G}_{\sigma _{i}}(f(\sigma _{i})\cap
 \mathfrak{F})=\mathfrak{G}_{\sigma _{i}}(h(\sigma _{i})\cap \mathfrak{F})\subseteq \frak{F}.$}

(2) {\sl $\mathfrak{F}=LF_{\sigma}(F)$, where $F$ is a
  formation $\sigma$-function
 such that 
$$F(\sigma _{i})=\mathfrak{G}_{\sigma _{i}}(f(\sigma _{i})\cap
 \mathfrak{F})=\mathfrak{G}_{\sigma _{i}}F(\sigma _{i})$$ for all $\sigma _{i}\in \Pi.$}

{\bf Proof.}  (1)  First suppose that   $\mathfrak{G}_{\sigma _{i}}(f(\sigma _{i})\cap
 \mathfrak{F}) \not \subseteq \mathfrak{F}$  and 
let $G$ be a group of minimal order in  $\mathfrak{G}_{\sigma _{i}}(f(\sigma _{i})\cap
 \mathfrak{F})  \setminus \mathfrak{F}$. Note that $f(\sigma _{i})\cap
 \mathfrak{F}$ is a formation  by Lemma  2.3(3), so 
 $\mathfrak{G}_{\sigma _{i}}(f(\sigma _{i})\cap
 \mathfrak{F})$ is a formation by  Lemma 2.1(1, 2). Hence 
 $R=G^{\mathfrak{F}}\leq  G^{f(\sigma _{i})\cap
 \mathfrak{F}}$  is a 
 unique minimal normal subgroup of $G$, so $R$   is a $\sigma _{i}$-group.

Moreover, $F_{\sigma _{i}}(G)=O_{\sigma _{i}}(G)$  and 
 $F_{\sigma _{j}}(G/R)=F_{\sigma _{j}}(G)/R $ for all $j\ne i$.  Therefore, since 
$G/R\in   \mathfrak{F}$ we have  $$(G/R)/F_{\sigma _{j}}(G/R)\simeq 
G/F_{\sigma _{j}}(G)\in   f(\sigma _{j})$$ for
 all $\sigma _{j}\in \sigma (G)\setminus \{\sigma _{i}\}$. Finally, we 
have   $$G/F_{\sigma _{i}}(G)=G/O_{\sigma _{i}}(G)  \in f(\sigma 
_{i})$$    
since  $G\in 
 \mathfrak{G}_{\sigma _{i}}(f(\sigma _{i})\cap
 \mathfrak{F})$ and the class ${\mathfrak{G}}_{\sigma _{i}}$ is hereditary. But then $G\in \mathfrak{F}$,
 a contradiction. Hence 
$\mathfrak{G}_{\sigma _{i}}(f(\sigma _{i})\cap \mathfrak{F})\subseteq 
\frak{F}.$

Now suppose that  $\mathfrak{G}_{\sigma _{i}}(f(\sigma _{i})\cap
 \mathfrak{F}) \not \subseteq \mathfrak{G}_{\sigma _{i}}(h(\sigma 
_{i})\cap \mathfrak{F})$  and 
let $G$ be a group of minimal order in   
 $\mathfrak{G}_{\sigma _{i}}(f(\sigma _{i})\cap
 \mathfrak{F}) \setminus \mathfrak{G}_{\sigma _{i}}(h(\sigma 
_{i})\cap \mathfrak{F}).$   Then   $G$ has a unique  minimal normal subgroup 
$R$,  $R=G^{\mathfrak{G}_{\sigma _{i}}(h(\sigma 
_{i})\cap \mathfrak{F})}$   and  $R\nleq O_{\sigma 
_{i}}(G)$. Hence $O_{\sigma _{i}}(G)=1$.

   Let $A$ be any non-identity $\sigma _{i}$-group and
 let 
$E=A\wr G=K\rtimes G$ be the regular wreath product of $A$ and    $G$,
 where $K$ is the base group of $E$. Then $O_{\sigma 
_{i}'}(E)=1$, so $F_{\sigma _{i}}(E)=O_{\sigma _{i}}(E)=K(O_{\sigma 
_{i}}(E)\cap G)=K$ since $O_{\sigma _{i}}(G)=1$. Moreover, since
 $G\in \mathfrak{G}_{\sigma _{i}}(f(\sigma _{i})\cap
 \mathfrak{F})\subseteq \mathfrak{F}$ we have  $E\in  \mathfrak{F}$ and so
 $E/F_{\sigma _{i}}(E)=E/K\simeq G\in h(\sigma _{i})\cap \mathfrak{F}
\subseteq 
\mathfrak{G}_{\sigma _{i}}(h(\sigma 
_{i})\cap \mathfrak{F}).$   Thus $\mathfrak{G}_{\sigma _{i}}(f(\sigma _{i})\cap
 \mathfrak{F})  \subseteq \mathfrak{G}_{\sigma _{i}}(h(\sigma 
_{i})\cap \mathfrak{F}),$ so $\mathfrak{G}_{\sigma _{i}}(f(\sigma _{i})\cap
 \mathfrak{F}) = \mathfrak{G}_{\sigma _{i}}(h(\sigma 
_{i})\cap \mathfrak{F}).$

(2)    Let $\mathfrak{M}=LF_{\sigma}(F)$.  Then   $$\mathfrak{M}=(\bigcap _{\sigma _{i}\in 
\Pi}{\mathfrak{G}}_{\sigma_{i}'} \mathfrak{G}_{\sigma _{i}}
({\mathfrak{G}}_{\sigma _{i}}(f(\sigma 
_{i})\cap \mathfrak{F}))) \cap  {\mathfrak{G}}_{\Pi }$$$$=(\bigcap _{\sigma _{i}\in 
\Pi}{\mathfrak{G}}_{\sigma_{i}'}
 {\mathfrak{G}}_{\sigma _{i}}(f(\sigma _{i})\cap \mathfrak{F})) \cap  {\mathfrak{G}}_{\Pi }=
\mathfrak{F}$$ by  Lemmas  2.3(3) and 2.4.
Hence we have (2).

The proposition    is proved. 

{\bf Corollary  2.6.}  (1) {\sl For every formation $\sigma$-function $f$ the 
class  $LF_{\sigma}(f)$ is a non-empty saturated formation.  }

(2) {\sl Every $\sigma$-local formation
    $\mathfrak{F}$  possesses a unique $\sigma$-local
 definition $F$ such 
that for every  $\sigma$-local definition $f$ of  $\mathfrak{F}$  and for 
every   $\sigma _{i}\in    \sigma ( \mathfrak{F})$   
 the following  holds:   }

{\sl $$F(\sigma _{i})=\mathfrak{G}_{\sigma _{i}}(f(\sigma _{i})\cap
 \mathfrak{F})=\mathfrak{G}_{\sigma _{i}}F(\sigma _{i}).$$  }

{\bf Proof.} (1)  First note that every identity group, by definition, 
belongs to  $LF_{\sigma}(f)$, so this class is  non-empty. On the other hand, 
the class  is a saturated  
formation by Lemma 2.3(3).

(2)  This assertion directly follows from Proposition 2.5(2).  

The corollary is proved.

Recall that $\text{form} (\mathfrak{X})$ denotes the intersection of all 
formations containing the collection of groups  $\mathfrak{X}.$

{\bf Proposition 2.7.}    {\sl Let $\mathfrak{F}=LF_{\sigma}(f)$ be a  
$\sigma$-local formation and $\Pi =\sigma (\mathfrak{F})$.  Let $m$ be the   
formation $\sigma$-function such that}  $m(\sigma _{i})=\text{form} 
(G/F_{\sigma _{i}}(G)   |\ G\in \mathfrak{F})$ {\sl for all  $\sigma _{i}\in \Pi$ and 
 $m(\sigma 
_{i})=\emptyset $ for all  $\sigma _{i}\in \Pi'$. Then:}

(i) {\sl $\mathfrak{F}=LF_{\sigma}(m)$, and }

(ii) {\sl $m(\sigma _{i}) \subseteq h(\sigma _{i})\cap \mathfrak{F}$ for every formation $\sigma$-function $h$
 of $\mathfrak{F}$ and for every $\sigma _{i}\in \sigma$.}

{\bf Proof. } Let
 ${\mathfrak{F}}(\sigma _{i})=(G/F_{\sigma _{i}}(G) |\ G\in 
{\mathfrak{F}})$  for all   $\sigma _{i}\in \Pi$,   
 and let $\mathfrak{M}=LF_{\sigma}(m)$. Then $\mathfrak{F}\subseteq 
\mathfrak{M}$. On the other hand, ${\mathfrak{F}}(\sigma _{i})\subseteq 
f(\sigma _{i})$  and  so $m(\sigma _{i})\subseteq f(\sigma _{i}) $ for all 
   $\sigma _{i}\in \Pi$. Also, we have $m(\sigma _{i})=\emptyset\subseteq f(\sigma _{i})$ for all
 $\sigma _{i}\in \Pi'$.  Hence  
$\mathfrak{M}\subseteq  \mathfrak{F}$, so $\mathfrak{M}=  \mathfrak{F}$.  
 The proposition   is proved. 

We call the $\sigma$-local definition $m$ of the formation  $\mathfrak{F}$ 
in Proposition 2.7 the \emph{smallest $\sigma$-local definition} of $\mathfrak{F}$.

\section{Proofs of Theorems 1.4 an 1.11}

 {\bf Proof of Theorem 1.11.}   (i) Suppose that  $\mathfrak{F}$  
  is  not  $\Sigma_{3}^{\sigma}$-closed and   
 let $G$ be a group of minimal order among the groups $G$ such that 
$G\not \in \mathfrak{F}$  but $G$ has subgroups $A_{1}$,   $A_{2}$ and 
$A_{3}\in  \mathfrak{F}$   such that
  the  indices 
$|G:A_{1}|$,  $|G:A_{2}|$ and   $|G:A_{3}|$
 are pairwise  $\sigma$-coprime.  Then $G=A_{i} A_{j}$ for all $i\ne j$.  
Let $R$ is a minimal normal subgroup of $G$.

(1) {\sl $G/R$ is $\sigma$-soluble $\Pi$-closed.  Hence $R$ is not a $\sigma$-primary
 $\Pi$-group.}

If for some $i$ we have $A_{i}\leq R$, then for any $j\ne i$ we have 
$G/R=A_{i}A_{j}/R=A_{j}R/R\simeq A_{j}/(A_{j}\cap R)\in 
\mathfrak{F}$ since $\mathfrak{F}$ is a formation by Lemma 2.2. Now assume that  $A_{i}\nleq R$ for all $i$.
Then  the hypothesis holds for $G/R$, so 
  $G/R$ is $\sigma$-soluble $\Pi$-closed by the choice 
of $G$.  Therefore, $R$ is  a not $\sigma$-primary
 $\Pi$-group since $G\not \in \mathfrak{F}$. Hence we have (1).

(2) {\sl $G$ is $\sigma$-soluble. }

Let  $L$ be a minimal normal subgroup of $A_{1}$.  Since $A_{1}$ is  $\sigma$-soluble,
  $L$ is a $\sigma _{i}$-group for some $i$. Moreover, since
 $|G:A_{2}|=|A_{1}:A_{1}\cap A_{2}|$ and $|G:A_{3}|=|A_{1}:A_{1}\cap A_{3}|$ are 
$\sigma$-coprime by hypothesis, we have   either $L\leq A_{1}\cap 
A_{2}$      or $L\leq A_{1}\cap 
A_{3}$. Therefore we can assume without loss of generality that $L\leq A_{2}$, so 
$L^{G}=L^{A_{1}A_{2}}=L^{A_{2}}\leq A_{2}$. Hence $1 < L^{G}$ is 
$\sigma$-soluble and so we have (2) by Claim (1).

(3) {\sl $R$ is a unique minimal normal subgroups of $G$, $R\nleq \Phi 
(G)$ and $R$ is a $\sigma _{i}$-group for some $\sigma _{i}\in \Pi'$.
 Hence $C_{G}(R)\leq R$.  }

Since $G$ is  $\sigma$-soluble by Claim (2),  $R$ is a $\sigma _{i}$-group for some 
$i$. Moreover, Claim (2) and Lemma 2.2 imply that
 $R$ is a unique minimal normal subgroups of $G$, $R$ is  a $\Pi'$-group  and
 $R\nleq \Phi (G)$.
 Hence $C_{G}(R)\leq R$ by \cite[A, 17.2]{DH}.

(4)  {\sl There are $j\ne k$ such that $R\leq A_{j}\cap  A_{k}$} (Since  
$|G:A_{j}|$ and  $|G:A_{k}|$   are  $\sigma$-coprime by hypothesis, this 
follows from Claim (3)).

{\sl Final contradiction for (i) } Since  $O_{\Pi}(A_{j})$ is 
  normal in   $A_{j}$ and 
$R\leq O_{\Pi'}(A_{j})$   by Claims  (3) and (4), we get that 
$O_{\Pi}(A_{j})\leq C_{G}(R)\leq R\leq O_{\Pi'}(A_{j})$ by Claim (3).
 Hence $  O_{\Pi}(A_{j})=1$. But $A_{j}$  is $\Pi$-closed by hypothesis and so   
$A_{j} $  is a  $\Pi'$-group. Similarly, one can show that
 $A_{k} $  is a  $\Pi'$-group and so $G=A_{j}A_{k}$ is a  
$\Pi'$-group. But then  $G$ is  is $\Pi$-closed. 
This contradiction  completes the proof of (i).

 (ii)  Suppose that $\mathfrak{M}$   
  is  not  $\Sigma_{3}^{\sigma}$-closed and   
 let $G$ be a group of minimal order among the groups $G$ such that 
$G\not \in \mathfrak{M}$  but $G$ has subgroups $A_{1}$,   $A_{2}$ and 
$A_{3}\in  \mathfrak{M}$   such that
  the  indices 
$|G:A_{1}|$,  $|G:A_{2}|$ and   $|G:A_{3}|$
 are pairwise  $\sigma$-coprime.     Then  $G\ne A_{i}$ for all $i$ 
 and   $G$ is $\sigma$-nilpotent by  Part 
(i). Moreover,  the choice 
of $G$  implies that 
  $G/R \in {\mathfrak{M}}$ for every minimal normal subgroup $R$ of $G$. 
Therefore  $R$ is a unique minimal normal subgroups of $G$ since  
$\mathfrak{M}$ is a formation. Hence, in 
fact,  $G$ is $\sigma _{i}$-group for some $i$. But then from $A_{1}  < G$ 
and $A_{2}  < G$  we get that   the  indices 
$|G:A_{1}|$ and  $|G:A_{2}|$ are not $\sigma$-coprime. This contradiction 
shows that we have (ii). The theorem is proved.

{\bf Lemma 3.1.}  {\sl If $G$ is $\sigma$-soluble,
 then $C_{G}(F_{\sigma}(G))\leq F_{\sigma}(G)$.}

{\bf Proof.} Let $C=C_{G}(F_{\sigma}(G))$. Assume that $C\nleq 
F_{\sigma}(G)$ and let $H/F_{\sigma}(G)$ be a chief factor of $G$ such that $H\leq
 F_{\sigma}(G)C$. Then  $H=F_{\sigma}(G)(H\cap C)$.     Since $G$ is $\sigma$-soluble,
 $$H/F_{\sigma}(G)=F_{\sigma}(G)(H\cap C)/F_{\sigma}(G)\simeq
 (H\cap C)/((H\cap C)\cap F_{\sigma}(G))$$  is a $\sigma _{i}$-group.

 Now 
let $U$ be a minimal supplement to  $(H\cap C)\cap F_{\sigma}(G)$ in  
$H\cap C$. Then  $((H\cap C)\cap F_{\sigma}(G))\cap  U\leq \Phi (U)$, so $U$ is
  $\sigma _{i}$-group. Moreover,  $(H\cap C)\cap F_{\sigma}(G)\leq Z(H\cap C)$ and  
so  $H\cap C$ is a 
normal ${\sigma}$-nilpotent subgroup of $G$. Hence $H\cap C\leq  F_{\sigma}(G)$ and so 
$H=F_{\sigma}(G)$.  This contradiction completes the proof of the lemma.

{\bf Lemma 3.2.}  {\sl  Let  $\mathfrak{F}= \mathfrak{S}_{\Pi}\mathfrak{X}$, where  
$\mathfrak{X}\subseteq \mathfrak{S}_{\sigma}$.   If the formation $\mathfrak{X}$ is 
$\Sigma_{t}^{\sigma}$-closed, then  $\mathfrak{F}$  is  $\Sigma_{t+1}^{\sigma}$-closed.}

{\bf Proof.}     Suppose that  
this lemma   is  false and   
 let $G$ be a group of minimal order among the groups $G$ such that 
$G\not \in \mathfrak{F}$  but $G$ has subgroups 
 $A_{1},  \ldots , A_{t+1}\in \mathfrak{F}$   such that 
  the indices
$|G:A_{1}|$,   \ldots ,  $|G:A_{t+1}|$  
 are pairwise  $\sigma$-coprime.   Then $G$ is $\sigma$-soluble by Theorem 
1.11. 

 Let $R$ be a minimal normal subgroup of $G$, so $R$ is a $\sigma _{i}$-group
 for some $i$. Moreover, the hypothesis holds for $G/R$ since $\mathfrak{F}$ is a formation by Lemma
 2.1(1, 2) and  so $G/R 
\in \mathfrak{F}$ by the choice of $G$. 
 Hence    
$R$ is a unique  minimal normal subgroup of                       
$G$ by the choice of $G$. Therefore, 
$\sigma _{i}\in \Pi'$ and $R\leq  O_{\sigma 
_{i}}(G)=F_{\sigma}(G)$.  Hence
 $C_{G}(F_{\sigma}(G)=C_{G}(O_{\sigma _{i}}(G))\leq O_{\sigma _{i}}(G)$  By Lemma 3.1.

  By hypothesis, there are numbers $i_{1},
 \ldots , i_{t}$ such
 that  $O_{\sigma _{i}}(G)\leq A_{i_{1}}\cap \cdots \cap  
A_{i_{t}}.$
Then $O_{\Pi}(A_{i_{j}})\leq C_{G}(O_{\sigma _{i}}(G))\leq O_{\sigma 
_{i}}(G).$
Hence $O_{\Pi}(A_{i_{j}})=1$ and so $A_{i_{j}}\in \mathfrak{X}$  for all $j=1, \ldots , t.$ 
Therefore $G\in \mathfrak{X}\subseteq \mathfrak{F}$ 
since  the formation $\mathfrak{X}$ is  $\Sigma_{t}^{\sigma}$-closed.
  This contradiction  completes the prove of the lemma.

  {\bf Lemma  3.3.}  {\sl Let  $\mathfrak{M}$ be a formation of $\sigma$-soluble   $\Pi$-closed 
groups and let  $\mathfrak{F}=\mathfrak{S}_{\Pi}\mathfrak{M}$.  If 
 $\mathfrak{M}$ is  $\Sigma_{3}^{\sigma}$-closed, then $\mathfrak{F}$  is
 $\Sigma_{3}^{\sigma}$-closed.}

{\bf Proof.}  Suppose that $G$ has subgroups 
  $A_{1}, A_{2},  A_{3}\in \mathfrak{F}$   such that  the  indices
$|G:A_{1}|$,   $|G:A_{2}|$, $|G:A_{3}|$
 are pairwise  $\sigma$-coprime. Then $G$ has a normal Hall 
$\Pi$-subgroup $V$  by Theorem 1.11.   Hence 
 $V\cap  A_{i}$ is a normal Hall $\Pi$-subgroup of  $A_{i}
$ and so from the isomorphism
 $VA_{i}/V\simeq A_{i}/A_{i}\cap V$ we get that $VA_{i}/V\in \mathfrak{M}$ and  
  the indices
$|(G/V):(A_{1}V/V)|$, $|(G/V):(A_{2}V/V)|$, $|(G/V):(A_{3}V/V)|$
 are pairwise  $\sigma$-coprime.  But then   $G/V\in 
\mathfrak{M}$ since $\mathfrak{M}$ is  $\Sigma_{3}^{\sigma}$-closed by hypothesis.  Hence   $G\in 
\mathfrak{F}$.  The lemma is proved.

The following lemma is evident.

{\bf Lemma 3.4.}  {\sl If  the class of groups $\mathfrak{F}_{j}$ is  
$\Sigma_{t}^{\sigma}$-closed for all $j\in J$, then  the class
 $\bigcap _{j\in J}\mathfrak{F}_{j}$  is also  $\Sigma_{t}^{\sigma}$-closed.}

A formation $\sigma$-function  $f$   is said to be: 
\emph{integrated} if $f(\sigma _{i})\subseteq LF_{\sigma}(f)$  for all $i$;  
 \emph{full} if $f(\sigma _{i})={\mathfrak{G}}_{\sigma _{i}}f(\sigma _{i})$   for all $i$.

In view of  Corollary  2.6, every  $\sigma$-local formation
    $\mathfrak{F}$  possesses a unique  integrated and full $\sigma$-local 
definition $F$.  We call such a  function $F$  the \emph{canonical $\sigma$-local definition}  
of  $\mathfrak{F}$.

{\bf Theorem 3.5.}  {\sl Let  $\mathfrak{F}=LF_{\sigma}(F)$ be a  $\sigma$-local formation of 
$\sigma$-soluble groups,  
 where $F$ is the canonical  $\sigma$-local definition    
of  $\mathfrak{F}$. If  the
 formation $F(\sigma _{i})$ is $\Sigma_{t}^{\sigma}$-closed for every $i$, then 
  the formation $\mathfrak{F}$ is $\Sigma_{t+1}^{\sigma}$-closed. }

{\bf Proof.} Let $\Pi =\text{Supp}(\mathfrak{F})$.  Then, by Lemma 2.3(4) 
and Corollary  2.6,   $$\mathfrak{F}= (\bigcap _{\sigma _{i}\in 
\Pi}\mathfrak{S}_{\sigma_{i}'} \mathfrak{G}_{\sigma _{i}}f(\sigma 
_{i})) \cap  \mathfrak{S}_{\Pi }= (\bigcap _{\sigma _{i}\in 
\Pi}\mathfrak{S}_{\sigma_{i}'} F(\sigma 
_{i})) \cap  \mathfrak{S}_{\Pi }.$$ 
  By Lemma 3.2, the formation $\mathfrak{S}_{\sigma _{i}'}F(\sigma 
_{i})$ is $\Sigma_{t+1}^{\sigma}$-closed. On the other hand, the class $\mathfrak{S}_{\Pi }$ is 
$\Sigma_{2}^{\sigma}$-closed  and so $\Sigma_{t+1}^{\sigma}$-closed.
  Hence   $\mathfrak{F}$ is 
 $\Sigma_{t+1}^{\sigma}$-closed  by Lemma 3.4. The  theorem is proved.

{\bf Proof of Theorem 1.4.}   Let  $\mathfrak{F}=LF_{\sigma}(f)$ be any 
 $\sigma$-local formation of 
meta-$\sigma$-nilpotent groups, where $f$ is the smallest  $\sigma$-local definition of 
$\mathfrak{F}$.  Then  the formation $f(\sigma _{i})$ is contained in
 ${\mathfrak{N}}_{\sigma}$ for all $\sigma _{i}
$ by Proposition 2.7.   Hence $f(\sigma _{i})$ is  $\Sigma_{3}^{\sigma}$-closed    by  Theorem 
1.11.   

Let $F$ be the canonical $\sigma$-local definition of  $\mathfrak{F}$.  
Then  $F(\sigma _{i})=\mathfrak{G}_{\sigma _{i}}f(\sigma _{i})$ for all $\sigma 
_{i}\in \sigma$ by Propositions 2.5 and 2.7.
  Hence, $F(\sigma _{i})$ is   
$\Sigma_{3}^{\sigma}$-closed  by Lemma 3.3.
 Therefore   $\mathfrak{F}$ is 
$\Sigma_{4}^{\sigma}$-closed  by  Theorem  3.5. 
   The theorem is proved.

\end{document}